\documentclass[12pt]{article}
\usepackage{amssymb}
\usepackage{theorem}

\newtheorem{theorem}{Theorem}[section]
\newtheorem{conjecture}[theorem]{Conjecture}
\newtheorem{proposition}[theorem]{Proposition}
\newtheorem{corollary}[theorem]{Corollary}
\newtheorem{lemma}[theorem]{Lemma}
{\theorembodyfont{\rmfamily}
\newtheorem{definition}[theorem]{Definition}

\newtheorem{problem}[theorem]{Problem}
\newtheorem{question}[theorem]{Question}
}

\include{epsf}

\begin{document}

\def\CP{\mathbb{C}{\rm P}}
\def\tr{{\rm tr}\,}
\def\endproof{{$\Box$}}

\title{Real line arrangements with Hirzebruch property}
\date{\today}
\author{Dmitri Panov \footnote{Supported by a Royal Society University Research Fellowship}}

\maketitle

\begin{abstract}

A line arrangement of $3n$ lines in $\mathbb CP^2$ satisfies Hirzebruch property if each line intersect others in $n+1$ points. Hirzebruch asked in \cite{Hir2} if all such arrangements are related to finite complex reflection groups. We give a positive answer to this question in the case when the line arrangement in $\mathbb CP^2$ is real, confirming that there exist exactly four such arrangements.
 
\end{abstract}

\section{Introduction and the main result }

The goal of this article is to prove the following result.

\begin{theorem}\label{main} There exist exactly four line arrangements  in $\mathbb RP^2$ consisting of $3\cdot n$ lines such that each line intersects others in $n+1$ points. These arrangements are reflection arrangements of the Coxeter groups corresponding to spherical triangles with angles $(\frac{\pi}{2},\frac{\pi}{2},\frac{\pi}{2} )$, $(\frac{\pi}{2},\frac{\pi}{3},\frac{\pi}{3} )$, $(\frac{\pi}{2},\frac{\pi}{3},\frac{\pi}{4} )$, $(\frac{\pi}{2},\frac{\pi}{3},\frac{\pi}{5})$.
\end{theorem}

Let us give a description of these four arrangements. The first arrangement is a union of there generic lines. The second arrangement is composed of three lines spanning the sides of a regular triangle in $\mathbb R^2$ together with three axes of symmetry of the triangle. The third arrangement is composed of four sides of a square in $\mathbb R^2$, four symmetry axes of the square, and the line at infinity. The fourth arrangement is composed of the sides of a regular pentagon in $\mathbb R^2$, five axes of symmetry, and five diagonals of the pentagon.

Following \cite{PP} we say that a line arrangement in $\mathbb CP^2$ satisfies {\it Hirzebruch property} if it consists of $3n$ lines and each line intersects others in exactly $n+1$ points. Such arrangements were studied first by Hirzebruch and H\"ofer in the context of construction of complex ball quotients\footnote{i.e. complex projective surfaces that are quotients of the unit complex ball $B^2_{\mathbb C}=\{|z_1|^2+|z_2|^2<1\}$ by a co-compact action of a discrete torsion free group.}. The ball quotients were obtained as disingularisations of ramified covers of $\mathbb CP^2$ with branching along line arrangements, the construction is described in \cite{Hir1} and \cite{BHH}.

Contemplating the list of arrangements suitable for construction of ball quotients Hirzebruch asked in \cite{Hir2} the following question:
\begin{question} Let  $\mathcal{L}$ be a complex line arrangement in $\CP^2$ consisting of $3\cdot n$ lines and such that each line of $\mathcal{L}$ intersect others at exactly $n+1$ points. Is it true that $\mathcal{L}$ is a {\it complex reflection arrangement}\footnote{A {\it complex reflection line arrangement} is a line arrangement in $\mathbb CP^2$ consisting of lines fixed by non-trivial elements of a finite complex reflection group acting on $\mathbb CP^2$. }?
\end{question}

This question  is still open, and Theorem \ref{main} gives a positive answer to it in the case when the line arrangement in $\mathbb CP^2$ is real.

Apart from the context of ball quotients, arrangements with Hirzebruch property appear in the setting of Polyhedral K\"ahler manifolds \cite{P}. This was used in  \cite{PP} to prove that the complement to any complex line arrangement with Hirzebruch property is aspherical. 

One more context in which these arrangements appear is the theory of  convex foliations on $\mathbb CP^2$, i.e. foliations whose leaves other than straight lines have no inflection points, see  Section \ref{discussion} and  \cite{MP} for more details.

{\bf About the proof.}  Theorem \ref{main} is deduced from existence of a special polyhedral metric with conical singularities on $\mathbb RP^2$ for which the lines of the arrangement are geodesics. The metric on $\mathbb RP^2$ is obtained by  restricting the polyhedral K\"ahler metric on the complexification of $\mathbb RP^2$, constructed in \cite{P} and whose properties are summarised in Section \ref{PKrecall}. To prove Theorem \ref{main} we show that the arrangement cuts $\mathbb RP^2$ into a collection of isometric  Euclidean triangles. Here we rely on a collection of elementary statements about spherical polygons, proven in Section \ref{equipolygons}. 

{\bf Acknowledgements.} I would like to thank Jorge Pereira and Anton Petrunin for useful and stimulating discussions and Piotr Pokora for comments on the first version of the paper.

\section{Polyhedral metrics}\label{polprop}

Recall the definition of polyhedral manifolds.
\begin{definition} Let $M$ be a  piecewise linear manifold $M$ with a complete metric $g$. We say that $M$ is a {\it polyhedral manifold of curvature} $\kappa$ if it admits a {\it compatible triangulation} for which each simplex equipped with $g$ is isometric to a geodesic simplex in the space of constant curvature $\kappa$. Depending on the sign of $\kappa$ the manifold $M$ is called a polyhedral spherical, Euclidean or hyperbolic manifold. The complement to metric singularities of a polyhedral manifold is denoted by $M^{\circ}$.
\end{definition}

Any polyhedral metric is non-singular in codimension $1$. The set of metric singularities  $M\setminus M^{\circ}$ is a union of some codimension two faces of a compatible triangulation. Let $\Delta$ be one of codimension two faces inside $M\setminus M^{\circ}$ and let $x$ be an interior point of $\Delta$.
Then in a neighbourhood of $x$ there is a geodesic surface orthogonal to $\Delta$ at $x$. The conical angle of such a surface at $x$ is the same for the all interior points of $\Delta$ and is called the {\it conical angle at} $\Delta$. 

We say that  a polyhedral Euclidean manifold $M$ is {\it non-negatively curved} if the conical angles at all  its codimesnion two faces are at most $2\pi$.

\subsection{Polyhedral surfaces}
A polyhedral surface is a polyhedral manifold of dimension two. Such a surface $S$ has a finite number of conical points $x_1,...,x_n$ and a complete metric $g$ which has constant curvature $\kappa$ on $S\setminus \{x_1,...,x_n\}$. In a neighbourhood of any conical point on $S$ there are polar coordinates $(r,\theta)$ in which  the metric can be given by the formulas $$g=dr^2+\alpha \kappa^{-1} \sin(\kappa r)d\theta^2,\;\; 
g=dr^2+\alpha rd\theta^2,\;\;g=dr^2-\alpha \kappa^{-1}\sinh(-\kappa r)d\theta^2,$$
depending on the sign of $\kappa$. The conical angle at $x$ is $2\pi\alpha$ in all these cases.

Each oriented polyhedral surface has a unique complex structure for which the polyhedral metric is K\"ahler on the complement to conical points. We will mainly study positively curved polyhedral metrics on $\mathbb CP^1$, invariant under the complex conjugation on $\mathbb CP^1$. Such metrics can be constructed by the {\it doubling} of {\it spherical polygons} that we will now describe.

{\bf Spherical polygons.} A {\it convex spherical polygon} is a closed convex subset of the sphere $\mathbb S^2_{\kappa}$ of curvature $\kappa$  with boundary composed of a finite number of geodesic segments. The geodesic segments are called the {\it edges} of the polygon and the points where these edges meet are called the {\it vertices}. If $P$ is a spherical (or Euclidean) polygon and $A$ is its vertex, we will denote the angle of $P$ at $A$ either by $\angle_A(P)$ or just by $\angle A$ (when the latter notation is unambiguous). We will assume that no two adjacent edges of the  polygon lie on one geodesic in $\mathbb S^2_{\kappa}$.

{\bf Doubling of polygons.} Let $P$ be a convex spherical polygon and let $P'$ be an isometric copy of it. The {\it doubling} of $P$ is obtained by gluing $P$ with $P'$ along their boundaries by the natural isometry. The resulting polyhedral sphere has a natural involution.

\begin{lemma} There is a one-to-one correspondence between convex spherical polygons and polyhedral metrics of positive curvature on $\mathbb CP^1$ satisfying the following properties.
\begin{itemize} 
\item The metric is invariant under the complex conjugation on $\mathbb CP^1$.
\item All the conical points are real, i.e., belong to $\mathbb RP^1\subset \mathbb CP^1$.
\item All the conical angles are less than $2\pi$.
\end{itemize}
\end{lemma}
 
The proof is straightforward, one direction of the correspondence is given by  the doubling construction. The other direction is given by taking the quotient of $\mathbb CP^1$ by the conjugation. Indeed, the conjugation is an isometry and so it leaves invariant a circle composed of geodesic segments.

\subsection{Polyhedral K\"ahler manifolds.}\label{PKrecall}
Here we recall  some definitions and results from \cite{P} concerning polyhedral K\"ahler manifolds.

\begin{definition}\label{PKdefinition} Let $\mathcal{M}$ be an orientable non-negatively curved Euclidean polyhedral manifold on dimension $2{\cdot}n$.
We say that  $\mathcal{M}$ is \emph{polyhedral K\"ahler} if the holonomy of the metric
on $\mathcal{M}^\circ$ belongs to $\mathrm{U}(n)\subset SO(2\cdot n)$.

\end{definition}

Our proof of Theorem \ref{main} relies heavily on the following theorem, proven in \cite{P}.

\begin{theorem}\label{hirmetric} Let $\cal L$ be an arrangement of $3n$ lines ($n\ge 2$) in $\mathbb CP^2$ with Hirzebruch property. Then there exists a unique up to scale polyhedral K\"ahler metric $g_{\cal L}^{\mathbb C}$ on $\mathbb CP^2$ which is singular along $\cal L$, non-singular in the complement of $\cal L$ and has conical angle $2\pi\cdot \frac{n-1}{n}$ at each line of the arrangement. 
\end{theorem}
The existence part of this theorem is a partial case of Theorem 1.12 in \cite{P}. The uniqueness of the metric up to scale follows from general results on unitary flat logarithmic connections.

{\bf The Euler field and the $S^1$-isometry.} It was proven in \cite{P}, that a polyhedral K\"ahler manifold $X$ of complex dimension two has the structure of a smooth complex surface, such that $X\setminus X^{\circ}$ is a divisor in $X$. Since $X$ is polyhedral, each point $x\in X$ has a conical $\varepsilon$-neighbourhood. It is obvious that on such a neighbourhood there is a real vector field $e_r$ acting by radial dilatation. In \cite{P} Section 3 it was explained that this field can be complexified to a holomorphic Euler field $e=e_r+ie_s$, and we sum up the properties of $e$ in the following theorem. It will be convenient to set $\varepsilon=2$ which can always be achieved by scaling the metric by a large factor.

\begin{theorem} Let $x\in X$ be a point, $U_x(2)$ be its conical neighbourhood of radius $2$, and $S_x(2)$ be the boundary of this neighbourhood. There is a holomoprhic Euler vector filed $e=e_r+ie_s$ defined on $U_x(2)$ with the following properties.

\begin{enumerate} 

\item The field $e_r$ is the real radial vector field acting by dilatations of the metric, it restricts to each ray of the cone as $r\frac{\partial}{\partial r}$.

\item The field $e_s$ is given by $e_s=J(e_r)$, where $J$ is  the  operator of complex structure on $TX$. The field $e_s$ is acting by isometries on $U_x(2)$.

\item  Let $x$ be a multiple point of an arrangement $\cal L$ from Theorem \ref{hirmetric} of \emph{multiplicity}\footnote{the multiplicity of a point is the number of lines of the arrangement passing through the point.} $\mu(x)\ge 2$. Then  $e_s$ integrates to an isometric $S^1$-action on $U_x(2)$ which is free on $U_x(2)\setminus x$. The quotient $S_x(2)/S^1$ is a curvature $1$ two-sphere with $\mu(x)$ conical singularities of angles $2\pi\cdot \frac{n-1}{n}$.

\end{enumerate} 

\end{theorem}

{\bf Proof.} This theorem is a partial case of Theorem 1.7 in \cite{P}.

\hfill $\square$

\subsubsection{Polyhedral K\"ahler metric for real line arrangements}

From now on we will assume that $\{L_1,...,L_{3n}\}=\cal L$ is a real line arrangement in $\mathbb RP^2$, satisfying  Hirzebruch property and $\{L_1^{\mathbb C},...,L_{3n}^{\mathbb C}\}=\cal L^{\mathbb C}$ is its complexification in $\mathbb CP^2$. Let $\sigma$ be the involution on $\mathbb CP^2$ induced by the complex conjugation, and let  $g_{\cal L}^{\mathbb C}$ be a polyhedral K\"ahler metric on $\mathbb CP^2$ given by Theorem \ref{hirmetric}, with conical singularities of angles $2\pi\frac{n-1}{n}$ at lines $L_i^{\mathbb C}$.

\begin{corollary}\label{realproperties} 

\begin{enumerate}

\item The polyhedral K\"ahler metric $g_{\cal L}^{\mathbb C}$ is invariant under the complex conjugation $\sigma$ on $\mathbb CP^2$. 

\item The metric $g_{\cal L}^{\mathbb C}$ restricts to a Euclidean polyhedral metric $g_{\cal L}^{\mathbb R}$ on $\mathbb RP^2$ and the lines $L_i$ are geodesics on $\mathbb RP^2$ with respect to $g_{\cal L}^{\mathbb R}$. 

\item Let $x$ be a real point $x\in \cal L\subset\cal L^{\mathbb C}$. Let $e=e_r+ie_s$ be the Euler field defined in a conical neighbourhood of $x$. Then $\sigma(e)=e_r-ie_s$.

\item The involution $\sigma$ descends to an isometry of the two-sphere $S_x(2)/S^1$, and $(S_x(2)/S^1)/\sigma$ is a convex spherical polygon of curvature $1$.
\end{enumerate}

\end{corollary}

{\bf Proof.}  1) The anti-holomorphic involution sends the polyhedral K\"ahler metric $g_{\cal L}^{\mathbb C}$ to a polyhedral K\"ahler metric. Since such a metric is unique up to scale, it is invariant under $\sigma$.

2) For any polyhedral metric the fixed set of any isometric involution is totally geodesic, so $\mathbb RP^2\subset \mathbb CP^2$ is totally geodesic. Hence the restriction of the metric to $\mathbb RP^2$ is a flat metric with conical singularities.  

To see that the lines $L_i$ are geodesic in $\mathbb RP^2$, note that each complex line $L_i^{\mathbb C}$ is totally geodesic in $\mathbb CP^2$, and $L_i$ is the fixed locus of the isometric involution $\sigma$ on $L_i^{\mathbb C}$. 

3) Let $e=e_r+ie_s$ be the holomorphic Euler field in a neighbourhood of $x$. Then $\sigma(e)$ is an anti-holomorphic vector field. A the same time, since $\sigma$ is an isometry preserving $x$, $\sigma(e_r)=e_r$. This proves the claim.

4) Indeed, from 3) it follows that $\sigma(e_s)=-e_s$, hence $\sigma$ sends $S^1$-orbits to $S^1$-orbits.

\hfill $\square$

\begin{definition}\label{polygonofapoint} For a real line arrangement $L_1,...,L_{3n}$ satisfying Hirzebruch property let $x$ be a multiple point. Denote by $\mathbb D(x)$ the convex spherical polygon  $(S_x(2)/S^1)/\sigma$ from  Corollary \ref{realproperties}.
\end{definition}

In the next Lemma we summarise what we need to know about polyhedral K\"ahler metrics in order to prove Theorem \ref{main}. 

Let ${\cal L}=\{L_1,...,L_{3n}\}$ be a real arrangement with Hirzebruch property. Suppose $x$ is a multiple point of ${\cal L}$ and assume that $k$ lines pass through $x$, i.e., $\mu(x)=k$. After a possible re-enumeration assume that the lines passing through $x$ are $L_1,...,L_k$ and they go in a cyclic order at $x$ on $\mathbb RP^2$. The spherical polygon $\mathbb D(x)$ associated to $x$ by Definition \ref{polygonofapoint} has $k$ vertices $A_1,...,A_k$ corresponding to the lines $L_1,...,L_k$. 

\begin{lemma}\label{mainlemma} The angle of the spherical polygon $\mathbb D(x)$ at each vertex $A_i$ is equal to $\pi\frac{n-1}{n}$. The angle between geodesics $L_i$ and $L_{i+1}$ on $\mathbb RP^2$ at the point $x$ with respect to the metric $g_{\cal L}^{\mathbb R}$ is equal to $\frac{1}{2}|A_iA_{i+1}|$ for all $i\in \{1,...,k\}$ (here $A_{k+1}=A_1$).
\end{lemma} 

{\bf Proof.} Let $U_x(2)$ be a conical $2$-neighbourhood of $x$ in $\mathbb CP^2$ with respect to the metric $g_{\cal L}^{\mathbb C}$. Consider its intersection with $\mathbb RP^2$, and let $S^1$ be the boundary of this intersection. Each line $L_i$ for $i\in \{1,...,k\}$ intersects $S^1$ in two points and we can denote them by $B_i$ and $B_{i+k}$, so that points $B_1,...,B_{2k}$ go along $S^1$ in a cyclic order. 

Denote by $\pi$  the quotient map $S_x(2)\to \mathbb D(x)$. Note that the map $\pi: S^1\to \partial(\mathbb D(x))$ is a locally isometric cover of degree two, and for any $i\in \{1,...,k-1\}$ the segment of $S^1$ included between $B_i$ and $B_{i+1}$ is sent isometrically to the edge $A_iA_{i+1}$ of $\mathbb D(x)$. Note finally that the length of $B_iB_{i+1}$ is twice the angle between $L_i$ and $L_{i+1}$ on $\mathbb RP^2$.   

\hfill $\square$

\section{Equiangular spherical  polygons} \label{equipolygons}

From now on by spherical polygons we mean polygons on the unit sphere $\mathbb S^2$. In the view of Lemma \ref{mainlemma} we will need to study equiangular spherical polygons.

\begin{definition}
A convex spherical polygon is called {\it equiangular} if the angles of the polygon at all vertices are equal. The polygon is called {\it equilateral} if all its edges are of the same length.
\end{definition}

The goal of this section is to prove the following proposition and its refinement Lemma \ref{uppebound} on equiangular spherical polygons.

\begin{proposition}\label{lessthan4pi3} Let $P^*$ be a  convex equiangular spherical polygon with $n\ge 3$ vertices. The sum of lengths of any two consecutive edges of $P$ is smaller than $\pi$ if $n$ is even and smaller than $2\pi-2\arccos(\frac{1}{n-1})$ if $n$ is odd.
\end{proposition} 

To each convex spherical polygon $P\subset \mathbb S^2$ with vertices $A_1,...,A_n$ one can associate the {\it dual convex polygon} $P^*$ with edges of lengths $\pi-\angle A_i$ and angles of value $|A_iA_{i+1}|$. To produce $P^*$ one  starts with the convex cone $C_P$ in $\mathbb R^3$ over $P\subset \mathbb S^2$, takes its dual cone $C_{P}^*$ and intersects it with $\mathbb S^2$, i.e., $P^*=C_P^*\cap \mathbb S^2$.  Clearly, this duality defines one-to-one correspondence between equiangular and equilateral polygons. So, Proposition \ref{lessthan4pi3} is equivalent to the following dual one, which we are going to prove.

\begin{proposition} \label{largerthan2pi3} Let $P$ be a  convex equilateral spherical polygon with $n\ge 3$ vertices. The sum of any two consecutive angles of $P$ is larger than $\pi$ if $n$ is even and greater than $2\arccos(\frac{1}{n-1})$ if $n$ is odd. 
\end{proposition} 

We will first reduce this statement to its Euclidean analogue by means of the following standard lemma.

\begin{lemma} For any convex spherical polygon $P$ with vertices $A_1,...,A_n$ there is a convex Euclidean polygon $P'$ with vertices $B_1,...,B_n$ such that for all $i$ $|A_iA_{i+1}|=|B_iB_{i+1}|$ and $\angle A_i>\angle B_i$.
\end{lemma} 

{\bf Proof.} Cut $P$ into $n-2$ convex triangles by diagonals $A_1A_i$. Replace each triangle by a flat one with the sides of the same length and glue back to get a flat polygon. Since the angles of all $n-2$ triangles have decreased, the resulting Euclidean polygon satisfies the condition of the lemma.

\hfill $\square$

To  prove Proposition \ref{largerthan2pi3} it remains to prove the following. 
\begin{proposition}\label{euclideanpol} Let $P$ be a  convex equilateral Euclidean polygon with $n\ge 3$ vertices. The sum of any two consecutive angles of $P$ is at least  $\pi$ if $n$ is even and at least $2\arccos(\frac{1}{n-1})$ if $n$ is odd. 
\end{proposition} 
This proposition in its turn will be deduced from the following two lemmas, the first of which is completely straightforward, and we omit its proof.

\begin{lemma}\label{deformation} For any convex Euclidean polygon $P$ with $n\ge 5$ vertices  $A_1,...,A_n$ there is an arbitrary small deformation of $P$ that preserves the lengths of edges and decreases the value $\angle A_1+\angle A_2$.
\end{lemma}

\begin{lemma}\label{quadrilateral} Let $ABCD$ be a convex Euclidean quadrilateral with sides of integer lengths such that $|AB|=1$ and $|AB|+|BC|+|CD|+|DA|=n$. Then $\angle A+\angle B\ge \pi$ if $n$ is even and  $\angle A+\angle B\ge 2\arccos(\frac{1}{n-1})$ if $n$ is odd.
\end{lemma}

{\bf Proof.} Consider first the case when $n$ is even. If $|CD|= 1$, $ABCD$ is a parallelogram, so we can assume  $|CD|>1$. There exists a unique parallelogram $ABC'D$ with $C'D=1$. Clearly,  
$\angle_A(ABC'D)=\angle_A(ABCD)$, and it is easy to check that $\angle_B(ABC'D)<\angle_B(ABCD).$
Since  $ABC'D$ is a parallelogram, we conclude $\angle_A(ABCD)+\angle_B(ABCD)>\pi$. 

Suppose now that $n$ is odd and assume $\angle A+\angle B<\pi$. Let $E$ be the intersection of the lines $\overline{AD}$ and $\overline{BC}$. Clearly 
$$|AC|+|CB|<|AD|+|DC|+|CB|=n-1<|AE|+|ED|,$$ 
so there is a point $F$ in the segment $EC$ such that $|AF|+|FB|=n-1$.
Clearly, $(\angle_A+\angle_B)(ABCD)>(\angle_A+\angle_B)(ABF)$. Note finally, that among all possible triangles of perimeter $n$ with one side of length $1$, the sum of two angles at this side attains its minimum for the isosceles triangle, and this minimum is $2\arccos(\frac{1}{n-1})$.

\hfill $\square$

{\bf Proof of Proposition \ref{euclideanpol}.} Let $\Pi_n$ be the space of all convex equilateral polygons in $\mathbb R^2$ 
with sides of length $1$. It has a natural compactification $\bar{\Pi}_n$ consisting of all convex polygons with sides of integer length. There is a continuous function $(\angle_{A_1}+\angle_{A_2})(P)$ defined on $\bar{\Pi}_n$ and from Lemma \ref{deformation} it follows that it attains its minimum on the part of $\bar{\Pi}_n$ consisting of quadrilaterals and triangles. Now the statement follows from Lemma \ref{quadrilateral}.

\hfill $\square$

The next lemma is a slight refinement of Proposition \ref{lessthan4pi3} for pentagons.

\begin{lemma}\label{uppebound} Any convex spherical equiangular pentagon satisfying $|A_{i-1}A_i|+|A_{i}A_{i+1}|>\frac{2\pi}{3}$ for $i=1,...,5$ satisfies $|A_{i-1}A_i|+|A_{i}A_{i+1}|<\pi$.

Dually, any convex spherical equilateral pentagon satisfying $\angle A_i+\angle A_{i+1}<\frac{4\pi}{3}$ for $i=1,...,5$ satisfies $\angle A_i+\angle A_{i+1}>\pi$.
\end{lemma} 

{\bf Proof.} Let us prove the dual statement. 
We will assume $\angle A_1+\angle A_2\le\pi$, and deduce that $\angle A_5+\angle A_1+\angle A_2+\angle A_3>\frac{8\pi}{3}$, which contradicts the conditions of the lemma. 

Let us decompose the pentagon into the union of the triangle $A_5A_4A_3$ and the quadrilateral  $A_5A_1A_2A_3$. The condition $\angle A_5+\angle A_1\le\pi$ implies $|A_1A_2|>|A_3A_5|$. So $|A_4A_5|=|A_4A_3|>|A_3A_5|$ and in the triangle $A_5A_4A_3$ the sum of angles at vertices $A_5$ and $A_3$ exceeds $\frac{2\pi}{3}$. Adding to this value the sum of all angles of the quadrilateral $A_5A_1A_2A_3$, which exceeds $2\pi$, we get the contradiction.

\hfill $\square$

The next lemma is straightforward, we omit the proof.
\begin{lemma}\label{equalregular} Let $k$ and $n$ be two integers with $n,k\ge 2$. Let $P_k$ be a regular (i.e, equilateral and equiangular) spherical $k$-gone and $P_n$ be a regular spherical $n$-gone. Suppose that the angles and the sides of $P_k$ have the same size as that of $P_n$. Then $n=k$.
\end{lemma}

\section{Proof of Theorem \ref{main}} \label{proof}

\subsection{Properties of the polyhedral metric $g_{\cal L}^{\mathbb R} $ on $\mathbb RP^2$}

Let us start the section by summarising the properties of the metric  $g_{\cal L}^{\mathbb R} $ on $\mathbb RP^2$ induced from the polyhedral K\"ahler metric  $g_{\cal L}^{\mathbb C} $ in $\mathbb CP^2$. First, we introduce some terminology. A real line arrangement $\cal L$ cuts $\mathbb RP^2$ into a collection of {\it polygons} whose edges are called {\it edges} of the arrangement. Two multiple points of $\cal L$ are called {\it adjacent} if they are the end points of one edge. 

For each multiple point $x$ of  $\cal L $ by the {\it star} $S(x)$ of $x$ we denote the union of all polygons adjacent to $x$. 
%Two lines are {\it adjacent} at $x$ if they contains two edges of one polygon in the star of $x$. 
The intersection of a small neighbourhood of $x$ with a star of $x$ is a union of $2\mu(x)$ {\it sectors}. 

\begin{theorem}\label{manyproperties} Consider a real line arrangement $\cal L$ of $3n$ lines with Hirzebruch property and let $g_{\cal L}^{\mathbb R}$ be the corresponding  metric on $\mathbb RP^2$. Then the following properties hold.

\begin{enumerate}

\item At any multiple point of $\cal L$ each sector has an acute angle unless the point is double, in which case all four sectors have angle $\frac{\pi}{2}$.

\item There is a constant $a(n)<\frac{\pi}{3}$ such that
the angles of sectors of all triple points of $\cal L$ are equal to $a(n)$. 

\item $\cal L$ is simplicial\footnote{i.e., all the polygons of the decomposition are triangles.}, and no two vertices of multiplicity two are adjacent.

\item Let $x$ be a multiple point of $\cal L$. The sum of angles of any two adjacent sectors of $x$ is less than $\frac{2\pi}{3}$ if $\mu(x)=3$, and less than $\frac{\pi}{2}$  if $\mu=4,5$.

\item  The multiplicity of each multiple point of $\cal L$ is at most $5$, and any point of multiplicity $5$ has exactly $5$ double points in the boundary of its star. 

\item For any multiple point of $\cal L$ the number of adjacent multiple points of multiplicity grater than $2$ is at most five.

\end{enumerate}
\end{theorem}

{\bf Proof.} Let $x$ be a multiple point of $\cal L$ and let $\mathbb D(x)$ be the associated spherical polygon. Its equiangular by Lemma \ref{mainlemma}.

1) The length of any edge of a convex spherical polygon is at most $\pi$ and it is equal to $\pi$ only in the case when the polygon is a bigon. Hence by Lemma \ref{mainlemma} the angle of each sector is at most $\frac{\pi}{2}$ and it is equal to $\frac{\pi}{2}$ iff $\mathbb D(x)$ has exactly two vertices, i.e., $x$ is a double point.

2) If $x$ is a triple point then $\mathbb D(x)$ is the unique spherical regular triangle with angles $\frac{\pi(n-1)}{n}$. The edges of such a triangle are shorter than $\frac{2\pi}{3}$, hence the statement holds by Lemma \ref{mainlemma}.

3) Since by property 1) the angles of all polygons in which the arrangement cuts $\mathbb RP^2$ are not obtuse, the only polygons different from triangles that can be present in the decomposition are rectangles. Assume by contradiction, that there is such a rectangle $R$ in the decomposition. Applying again property 1) we see that all vertices of $R$ are double points. If follows that all polygons sharing an edge with $R$ are rectangles as well. Applying this reasoning repeatedly we come to contradiction.

4) This is a direct consequence of Proposition \ref{lessthan4pi3} and Lemma \ref{uppebound} applied to the polygon $\mathbb D(x)$.

5) Let $x$ be a point of the arrangement of multiplicity $d$ and let $S(x)$ be its star. This star is a union of triangles by property 3). Denote by $P_1,P_2,...,P_{2d}$ the vertices of these triangles lying on the boundary of $S(x)$, enumerated in a cyclic order. Note that unless the point $P_i$ is a double point of the arrangement, by property 4) the angle of $S(x)$ at $P_i$ is less than $\frac{2\pi}{3}$. We deduce from 3) that there are at least $d$ points in the boundary of $S(x)$ with angle less than $\frac{2\pi}{3}$. Since the boundary of $S(x)$ is convex and the conical angle at $x$ is less than $2\pi$, applying Gauss-Bonnet formula to the star $S(x)$ we conclude that $d\le 5$.

6) The proof of this statement is identical to the proof of statement 5).

\hfill $\square$

\subsection{Proof of Theorem \ref{main}}

To prove Theorem \ref{main} we will show that all the triangles in the decomposition of $\mathbb RP^2$ by $\cal L$ are isometric with respect to the metric $g_{\cal L}^{\mathbb R}$. We will start with the following lemma.

\begin{lemma}\label{no4455} Let $x$ and $y$ be two adjacent multiple points in a real arrangement satisfying Hirzebruch property. Suppose $\mu(x),\mu(y)\ge 3$. Then $\mu(x)=3$ or $\mu(y)=3$.
\end{lemma}

{\bf Proof.} Consider triangles $\Delta_1$ and $\Delta_2$ of the decomposition that contain the edge $xy$ and let $Q_1$ and $Q_2$ be their vertices opposite to $xy$. Since the angles at points $Q_1$ and $Q_2$ can not be obtuse by Theorem \ref{manyproperties} 1), in quadrilateral $xQ_1yQ_2$ we have: $\angle x+\angle y\ge \pi$. Hence either $\angle x\ge \frac{\pi}{2}$ or $\angle y\ge \frac{\pi}{2}$, and the corresponding point is of multiplicity three by Theorem \ref{manyproperties} 4-5).

\hfill $\square$

The next two corollaries give a complete description of stars vertices of multiplicities $4$ and $5$.

\begin{corollary}\label{fivethree} Let $x$ be a point of multiplicity five of a real arrangement with Hirzebruch property. Let $P_1,...,P_{10}$ be the multiple points of the arrangement at the boundary of $S(x)$ and assume that $\mu(P_1)=2$. Then for $i=1,...,5$ we have $\mu(P_{2i-1})=2$ and $\mu(P_{2i})=3$.  
\end{corollary}

{\bf Proof.} By Theorem \ref{manyproperties}  5) five of points $P_1,...,P_10$ have multiplicity $2$. Hence it follows Theorem \ref{manyproperties} 3) that points $P_{2i-1}$ have multiplicity two. The remaining five points have multiplicity $3$ by Lemma \ref{no4455}.

\hfill $\square$

\begin{corollary}\label{fourpoint} Suppose $x$ is a point of multiplicity four of a real arrangement with Hirzebruch property, and let $P_1,...,P_8$ be the vertices of its star. Then at least one points $P_i$, say $P_1$, has multiplicity $2$. In such a case for $i=1,...,4$ we have $\mu(P_{2i-1})=2$, $\mu(P_{2i})=3$.
\end{corollary}

{\bf Proof.} To see that $x$ has at least one adjacent point of multiplicity $2$ we apply Theorem \ref{manyproperties} 6). Denote this double point by $P_1$.

By Lemma \ref{no4455} points $P_1,...,P_8$ can not have multiplicity four of five. So it is enough to show that there can not be five points of multiplicity $3$ in the star of $x$. Assume the contrary. In this case there are three consequent vertices, say $P_2, P_3, P_4$ of multiplicity three. Consider triangles $xP_2P_3$ and $xP_3P_4$ and note that the angles at points $P_i$ are less than $\frac{\pi}{3}$ by Theorem \ref{manyproperties} 2). Hence by Gauss-Bonnet $\angle_x(xP_2P_3)+\angle_x(xP_3P_4)>\frac{2\pi}{3}$, which is a contradiction with Theorem \ref{manyproperties} 4).

\hfill $\square$

An immediate consequence of Corollaries \ref{fivethree} and \ref{fourpoint} is the following statement.
\begin{corollary}\label{equalangles} Let $\cal L$ be a real line arrangement with Hirzebruch property and let $x$ be its multiple point. All sectors at $x$ have the same angle at $x$ with respect to the metric $g_{\cal L}^{\mathbb R}$.
\end{corollary}

{\bf Proof.} If $x$ is a double or triple point then this statement holds by Theorem \ref{manyproperties}. 

Suppose $x$ is a point of multiplicity $4$. Using notations of Corollary \ref{fourpoint} we see that for any $i=1,...,7$ triangles $xP_iP_{i+1}$ and $xP_{i+1}P_{i+2}$ ($P_9=P_1$) are isomteric by an isometry that sends  $P_i$ to $P_{i+2}$ and  fixes $P_{i+1}$  and $x$. Hence all $8$ sectors at $x$ have the same angle. 

The case $\mu(x)=5$ follows from Corollary \ref{fivethree} in the same way.

\hfill $\square$

\begin{corollary}\label{alltriangles} Suppose that $x$, $y$, $z$ are adjacent points of a real arrangement with Hirzebruch property. Then the multiplicities of these points belong the the following list (up to a permutation): $(2,3,3)$, $(2,3,4)$, $(2,3,5)$.
\end{corollary}
{\bf Proof.} By Lemma \ref{no4455} at most one of points $x$, $y$, $z$ can have multiplicity $4$ or $5$. Assume that this point is $z$. Then applying to the star of $z$ either Corollary \ref{fivethree} or Corollary \ref{fourpoint} we see that multiplicities of $x$ and $y$ are $(2,3)$ up to a permutation. 

All three points of the triangle $xyz$ can't be of multiplicity $3$ since in this case $\angle x=\angle y=\angle z<\frac{\pi}{3}$ by \ref{manyproperties} 2), which contradicts Gauss-Bonnet.

\hfill $\square$

\begin{corollary}\label{isotriangles} Let $\cal L$ be a real line arrangement with Hirzebruch property.

 \begin{enumerate}

\item The lines of $\cal L$ cut $\mathbb RP^2$ into isometric triangles with respect to the metric $g_{\cal L}^{\mathbb R}$. 

\item There is some $d\in \{3,4,5\}$ such that the multiplicities of vertices of each triangle are $(2,3,d)$ up to a permutation.
 
\end{enumerate}
\end{corollary}
{\bf Proof.} 1) Let $xyz$ and $xyt$ be two triangles of the decomposition that share the side $xy$. Then by Corollary \ref{equalangles} these triangles have the same angles at $x$ and $y$. Hence they are isometric. Hence all triangles of the decomposition are isometric.
% 1) Let $xyz$ and $xyt$ be two triangles of the decomposition that share the side $xy$. Suppose first that one of points $x$ or $y$, say $x$, has multiplicity $4$ or $5$. In this case by Corollary \ref{alltriangles} points $z$ and $t$ have the same multiplicity equal to two or tree. Then by Theorem \ref{manyproperties} 1-2) we have $\angle_z(xyz)=\angle_t(xyt)$. At the same time the point $y$ has multiplicity $2$ or $3$, and so $\angle_y(xyz)=\angle_y(xyt)$. Hence the triangles $xyz$ and $xyt$ are isometric. 

%The remaining case is when $x$ and $y$ have multiplicities $2$ and $3$ (or $3$ and $2$) correspondingly. Then $\angle_x(xyz)=\angle_x(xyt)=\frac{\pi}{2}$, and $\angle_y(xyz)=\angle_y(xyt)$. Hence again $xyz$ and $xyt$ are isometric.

2) By Corollary \ref{alltriangles} for any two triangles of the decomposition there vertices can be denoted by $x,y,z$ and $x',y',z'$ so that
 $$\mu(x)=\mu(x')=2,\;\;\mu(y)=\mu(y')=3,\;\;\mu(z)=d,\;\;\mu(z')=d', \;\; d,d'\ge 3.$$
In this case by 1) there is an isometry between the triangles that sends $x$ to $x'$, $y$ to $y'$ and $z$ to $z'$. By Corollary \ref{equalangles} the spherical polygons $\mathbb D(x)$ and $\mathbb D(x')$ are regular. Moreover, since $\angle_{z}=\angle_{z'}$, the polygons have sides of same the length and additionally they have angles of size $\frac{(n-1)\pi}{n}$ by Lemma \ref{mainlemma}. Hence $d=d'$ by Lemma \ref{equalregular}.

\hfill $\square$

{\bf Proof of Theorem \ref{main}.} According to  Corollary \ref{isotriangles} we have $3$ cases $d=3,4,5$. Replace each triangle in $\mathbb RP^2$ by 
a spherical triangle (of curvature $1$) with angles $(\frac{\pi}{2},\frac{\pi}{3},\frac{\pi}{d})$. As a result we obtain an $\mathbb RP^2$ with curvature $1$ metric and a Coxeter arrangement in it.

\hfill $\square$

\section{Discussion}\label{discussion}

In article \cite{Hir2} Hirzebruch gives the list of complex reflection arrangements of $3n$ lines, such that each line intersect others in $n+1$ points. This list consists of two infinite series and five exceptional examples. The infinite series are called $A_m^0$  or Ceva arrangements ($m\ge 3$) and $A_m^3$ $(m\ge 2)$ (or extended Ceva arrangements) and correspond to reflection groups $G(m,m,3)$ and $G(m,p,3)$ ($p<m$) from Shephard-Todd classification.
The arrangements $A_m^0$ and $A_m^3$ are given in homogeneous coordinates by equations 
$$(z_0^m-z_1^m)(z_1^m-z_2^m)(z_2^m-z_0^m)=0,$$
$$z_0z_1z_2(z_0^m-z_1^m)(z_1^m-z_2^m)(z_2^m-z_0^m)=0\;\;\;$$
correspondingly. The five  exceptional examples are associated to reflection groups $G_{23} , G_{24}, G_{25}, G_{26}, G_{27}$.  The corresponding arrangements are called the icosahedron configuration ($15$ lines),  the configuration $G_{168}$ or Klein configuration ($21$ lines), the Hesse configuration ($12$ lines), the configuration $G_{216}$ or extended Hesse configuration ($21$ lines),  and the configuration $G_{360}$ or Valentiner configuration  ($45$ lines), see \cite{Hir2}.

I believe that in the view of Theorem \ref{main} one can restate Hirzebruch's question as a conjecture.
\begin{conjecture}All arrangements satisfying Hirzebruch property are complex reflection arrangements.
\end{conjecture}

{\bf Convex foliations.} Line arrangements with Hirzebruch property have an interesting relation to {\it reduced convex foliations} in $\mathbb CP^2$. A foliation  in $\mathbb CP^2$ is called {\it convex} is its leaves other than straight lines have no inflection points. A foliation is called {\it reduced} if its inflection divisor is reduced \cite{MP}. It turns out, that any arrangement which can be realised as the union of all lines tangent to a reduced convex foliation, satisfies Hirzebruch property. Moreover all arrangements from Hirzebruch's list apart from $G_{169}$ and $A_6$ are indeed realised as line arrangements of reduced convex foliations (see \cite{MP} for more details).

It was explained \cite{Per} that any real line arrangement realisable as the line arrangement of a convex foliation is simplicial, which can be seen as a partial case of Theorem \ref{manyproperties} 3).
Note that at the present only a conjectural classification of simplicial arrangements in $\mathbb RP^2$ is known, see \cite{G1}, \cite{G2}.

{\bf Real polyhedral K\"ahler metrics.} Theorem \ref{main} can be seen as a first step toward a solution of the following classification problem. 

\begin{definition} A polyhedral K\"ahler metric on $\mathbb CP^2$ is called real if it is invariant under the conjugation of $\mathbb CP^2$. We call this metric {\it maximally real} if  the divisor of singularities of the metric is smooth in the complement of $\mathbb RP^2$.
\end{definition}

\begin{problem} Classify all positively curved maximally real polyhedral K\"ahler metrics on $\mathbb CP^2$.
\end{problem}


\begin{thebibliography}{99999999999999}


\bibitem[BHH]{BHH} G. Barthel, F. Hirzebruch, T. H\"ofer. Geradenkonfigurationen und Algebraische Fl\"achen.  Aspects of Mathematics, D4. Friedr. Vieweg \& Sohn, Braunschweig, (1987).

\bibitem[G1]{G1} B. Gr\"unbaum. Arrangements of hyperplanes.
Proceedings of the second Louisisana conference on combinatorics 
graph theory and computings. 41--106 (1971).

\bibitem[G2]{G2} B. Gr\"unbaum. Arrangements and spreads. Regional
Conf. Series in Mathematics {\bf 10}. Amer. Math Soc., (1972).   

\bibitem[Hir1]{Hir1} F. Hirzebruch. Arrangements of lines and algebraic surfaces. Arithmetic and Geometry. Volume 36 of the series Progress in Mathematics (1983), 113--140.

\bibitem[Hir2]{Hir2} F. Hirzebruch. Algebraic surfaces with extreme Chern numbers.
{\it Russian Math. Surveys} {\bf 40} (1985), 135--145.


\bibitem[MP]{MP}  David Marin, J. V. Pereira. Rigid flat webs, Asian Journal of  Mathematics 17 (2013), no. 1, 163--192.

\bibitem[P]{P} D. Panov. Polyhedral K\"ahler manifolds. Geometry and Topology 13 (2009) 2205--2252.

\bibitem[PP]{PP} D. Panov, A. Petrunin. Ramification conjecture and Hirzebruch's property of line arrangements arXiv:1312.6856. To appear in Compositio Mathematica.

\bibitem[Per]{Per} J. V. Pereira. Convex foliations, perprint.


%\bibitem[Th]{Th} W. P. Thurston. Shapes of polyhedra and triangulations of the sphere.{\it Geometry and Topology Monographs} {\bf 1}, (1998) 511--549.

%\bibitem[Tr]{Tr} M. Troyanov. Les surfaces euclidiennes {\`a}singularit{\'e}s coniques. {\it L'Enseign. Math.} {\bf 32}(1986), 79--94.

\end{thebibliography}
\end{document}